\documentclass{article}

\usepackage{arxiv}
\usepackage[utf8]{inputenc} 
\usepackage[T1]{fontenc}    
\usepackage{hyperref}       
\usepackage{url}            
\usepackage{booktabs}       
\usepackage{amsfonts}       
\usepackage{nicefrac}       
\usepackage{microtype}      
\usepackage{graphicx}
\usepackage{natbib}
\usepackage{doi}
\usepackage{nicefrac}
\usepackage{mathtools,amsmath,amsfonts,enumerate,nicefrac,graphicx,caption,subcaption,amssymb, xcolor, array}
\usepackage{comment}

\usepackage{titlesec}
\titlespacing{\section}{0pt}{12pt}{0pt}
\titlespacing{\subsection}{0pt}{8pt}{0pt}
\titlespacing{\subsubsection}{0pt}{8pt}{0pt}

\usepackage[page]{appendix}

\usepackage{newtxtext,newtxmath}
\renewcommand{\arraystretch}{1.2}

\usepackage[affil-it]{authblk} 
\usepackage{lipsum} 

\usepackage{fancyhdr}
\usepackage{makecell}

\usepackage{multirow}
\usepackage{makecell}
\renewcommand{\arraystretch}{1.5}

\usepackage{natbib}
\usepackage{verbatim}

\usepackage[nameinlink]{cleveref}
\usepackage{float}

\title{Can the perfect Swiss alphorn be designed? \\
A combination of reduced basis method and machine learning for shape optimization}

\author[1,\thanks{Corresponding author's email: fabio.marcinno@epfl.ch}]{Fabio Marcinn\`o}
\author[1]{Clément Froidevaux}
\author[1]{Leonardo Bocchieri}
\author[1]{Riccardo Tenderini}
\author[2]{Gerald Pot}
\author[1]{Simone Deparis}

\affil[1]{Institute of Mathematics, École Polytechnique Fédérale de Lausanne (EPFL), Lausanne, Switzerland}
\affil[2]{Alphorn Atelier, Monthey, Switzerland}
 
\date{}
\renewcommand{\undertitle}{}

\hypersetup{
pdftitle={Can the perfect Swiss alphorn be designed?},
pdfsubject={Wind Instrument, Helmholtz Equation, Shape Optimization, Reduced Basis Method, Machine Learning},
pdfauthor={Fabio Marcinnò, Clément Froidevaux, Leonardo Bocchieri, Riccardo Tenderini, Gerald Pot, Simone Deparis},
pdfkeywords={Wind Instrument, Helmholtz Equation, Shape Optimization, Reduced Basis Method, Machine Learning},
}
\begin{document}
\maketitle

\begin{abstract}
In this work, we investigate the shape optimization of a Swiss alphorn to achieve resonance frequencies as close as possible to prescribed target notes.
We construct an accurate geometric model of the alphorn based on both experimental measurements and imaging techniques. The latter are used to reconstruct geometric information where direct measures are not available. The resulting 3D model is meshed using a fully structured approach.
Then, we introduce a geometry-parametrized finite element (FE) formulation of the Helmholtz equation enabling the application of the Reduced Basis Method (RBM) to efficiently generate a large dataset of simulations, in which each set of geometric parameters corresponds to a specific set of resonance frequencies (notes). This dataset enables the training of machine learning (ML) models for forward prediction of resonance frequencies from geometric parameters, as well as inverse design, where geometries are estimated to achieve specified target notes.
\end{abstract}


\section{Introduction}\label{sec:intro}
The development of traditional acoustical music instruments is historically a slow process, made of continuous trials and errors. Every small innovation requires the work of many generations of makers who have to effectively fabricate prototypes to verify their findings. This is especially true for those kinds of instruments whose essence is inevitably linked to a local culture, like the Swiss alphorn \cite{kammermann2020}.
In this context, computational studies are extremely useful for accelerating the design process, testing acoustic hypotheses and gaining deeper insights into the physical behavior of the instrument without the need for physical prototyping, see \cite{fehlmann1994}, \cite{noreland2003}, \cite{lefebvre2011} and \cite{gerasimov2023}.

\subsection{The Swiss Alphorn}\label{sec:alphorn}
\begin{figure}[t]
        \centering
        \includegraphics[width=1\linewidth]{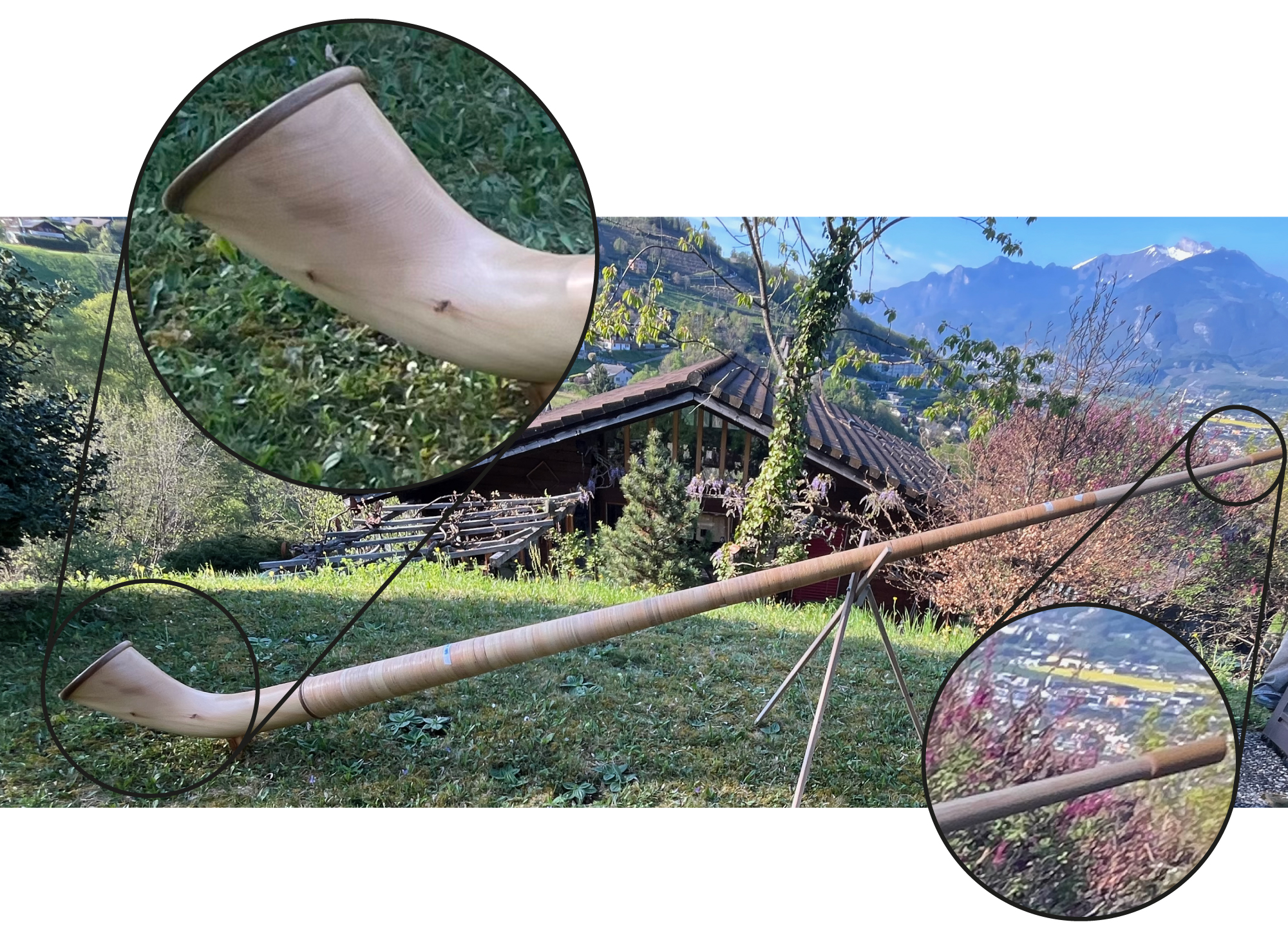}
        \caption{The Swiss alphorn along with a zoom on the bell and mouthpiece.}
        \label{fig:alphorn}
    \end{figure}

Our study is centered on the Swiss alphorn, a wind musical instrument with an overall conical shape, starting with a mouthpiece and ending with a curved bell as shown in Fig.\ \ref{fig:alphorn}. The instrument is made of softwood, in particular of spruce where air is blown inside the tube with the help of a cup-shaped, air-jet mouthpiece. Unlike reed instruments, sound is produced by the vibration of the lips directly against the wooden surface of the mouthpiece. The instrument has no toneholes on its surface, and therefore cannot play the entire chromatic scale of notes, i.e.\ the set of twelve semitones within an octave composing the traditional Western music. Since the design of this instrument has historically been a personal secret of each horn maker, there are still no standard specifications for the length or bore diameter, nor any imposed guidelines regarding its assembly \cite{kammermann2020}.

\subsection{Literature Review}\label{sec:literature}
Generally speaking, acoustic waves travel inside the tube of the instrument, from the tip to the first opening, partially reflecting back toward the mouthpiece. These vibrations generate standing waves, which are associated with specific musical notes. Any perturbation in the shape of the bore can affect the wave shape, and thus the frequency of the vibration.
The acoustic response of an instrument is characterized by its input
impedance evaluated at the mouthpiece. The impedance is a complex-valued, frequency-dependent quantity computed from the acoustic pressure field.
Resonances of the instrument appear as sharp peaks in the magnitude of the impedance.
These peaks correspond to frequencies at which the air column efficiently sustains
oscillations \cite{nederveen1998}. 
In the last decades, a lot of effort has been spent to model the input impedance and calculate its peaks, corresponding, as mentioned, to the natural frequencies at which the instrument vibrates the most.
Early important results for wind instruments, such as \cite{bakus1974} and \cite{causse1984}, focus on the experimental calculation of the input impedance.

A well-established approach to the modeling of wind instrument impedance is the Transmission-Matrix Method (TMM), discussed in \cite{plink1979}, \cite{keefe1990} \cite{bowen2018}. 
The TMM method approximates the geometry of the instrument as a sequence of concatenated segments, such as in the case of circuit components.
Due to its limitations in the modeling of complex geometries, this method has gradually been supplemented with FEM models, see \cite{noreland2003}, \cite{lefebvre2011}, \cite{gerasimov2023}.
Recently, to reduce the computational cost associated with FEM, methods such as Boundary Elements Methods (BEM) have been adopted for the numerical modeling of the sound propagation in ducts \cite{kreuzer2022}.
All the mentioned works deal with equations in the frequency domain, ideal for resonance analysis.

Other relevant methods, such as Finite-Difference Time-Domain (FDTD) used in \cite{redondo2004}, \cite{bilbao2013} and Lattice Boltzmann Method (LTB) adopted in \cite{gabriel2020}, are based on a time domain framework, particularly effective for analyzing transient acoustic behavior, which is not within the scope of this study.
Moreover, as explained in \cite{lefebvre2011}, the acoustic waves of wind instruments are supposed to be of sufficiently small amplitude for nonlinearities to be negligible, so we have decided to adopt the stationary and linear acoustic equation in the frequency domain, namely the Helmholtz equation.

Concerning the influence of the material, \cite{widholm2001} and \cite{nief2008} agree that the material and thickness of the instrument wall do not influence the main vibro-acoustics properties, conversely affecting mostly the musical timbre. There is an agreement on considering the acoustic properties of woodwind instruments mainly a consequence of its geometry.
For this reason, an analysis of the problem using FEM in combination with a geometry-based Reduced Basis Method (RBM) and machine learning represents a possible optimal way to approach the shape optimization of the Swiss alphorn.

\subsection{Novelties and Outline}\label{sec:novelties_and_outline}
The objective of this work is to integrate reduced-order modeling and machine learning techniques for the alphorn shape optimization, since finite element simulations, although accurate, are computationally expensive and thus impractical for direct use in such studies.
Firstly, we describe the reconstruction of the 3D Swiss alphorn model and the generation of its numerical mesh.
Secondly, a geometry-based parametrized FE formulation of the Helmholtz equation is derived on a reference domain, enabling the application of the RBM to efficiently generate a large dataset of simulations, where each set of geometric parameters corresponds to certain resonance frequencies (notes).
Finally, the dataset is used for training machine learning models addressing both forward prediction of resonance frequencies and inverse estimation of geometries producing prescribed target notes.

The work is organized as follows. 
In Sec.\ \ref{sec:preprocessing_and_mesh_generation}, we present the pre-processing needed for the reconstruction of the 3D model and the mesh generation.
Sec.\ \ref{sec:mathematical_modeling} presents the mathematical model along with the RBM framework and the geometry-based parametrization of the equation. Numerical settings and results are reported in Sec.\ \ref{sec:results}. Finally, conclusions and limitations are discussed in Sec.\ \ref{sec:conclusions}.

\section{Pre-Processing and Mesh Generation} \label{sec:preprocessing_and_mesh_generation}
\begin{figure}[t]
    \centering
    \includegraphics[width=1\linewidth]{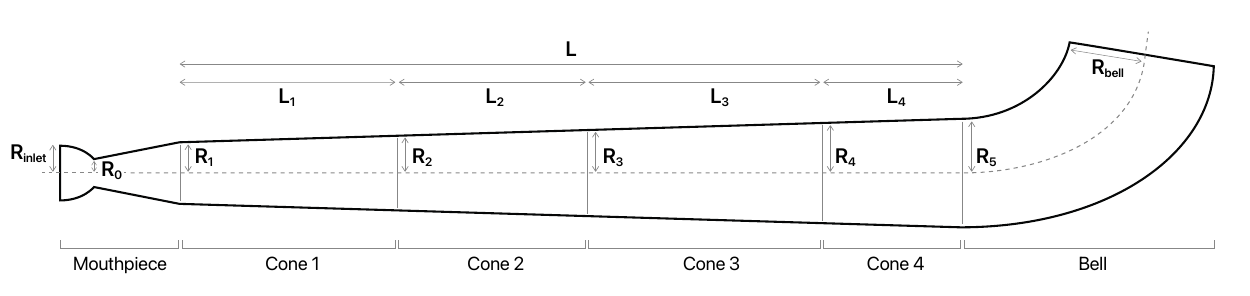}
    \caption{Subdivision of the alphorn in 6 principal pieces. $L_{k}$ are the lengths of the pieces and $R_{k}$ are the radii. The mouthpiece and the cones are axisymmetric.}
    \label{fig:alphorn_division}
\end{figure}

To accurately simulate the impedance spectrum and therefore the resonance frequencies of the instrument, the first step is to carefully reconstruct the 3D model.

\begin{figure}[t!]
    \centering
    \begin{subfigure}[t]{0.55\linewidth}
        \centering
        \includegraphics[width=\linewidth]{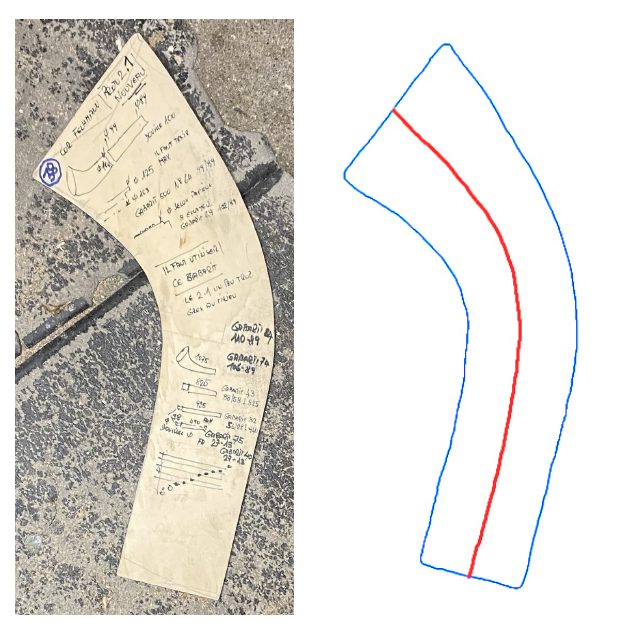}
        \caption{Mold of the bell and extraction of the contours and centerline.}
        \label{fig:pre-processing:steps1}
    \end{subfigure}
    \hfill
    \begin{subfigure}[t]{0.38\linewidth}
        \centering
        \includegraphics[width=\linewidth]{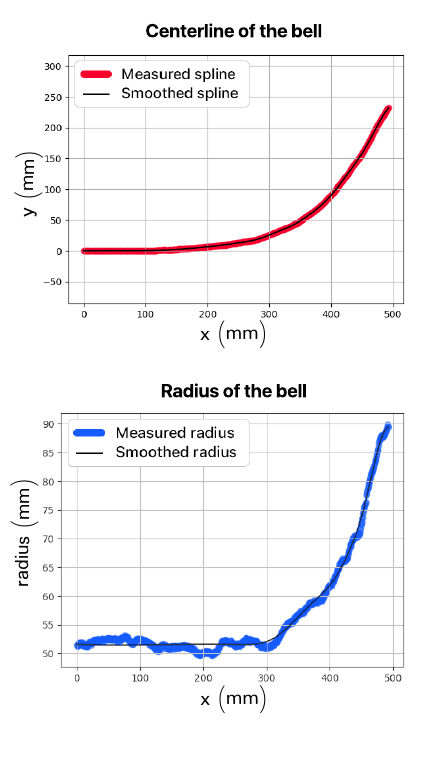}
        \caption{Smoothing of the geometry of the bell.}
        \label{fig:pre-processing:steps2}
    \end{subfigure}
    \vspace{0.5cm}

    \begin{subfigure}[t]{\linewidth}
        \centering
        \includegraphics[width=0.75\linewidth]{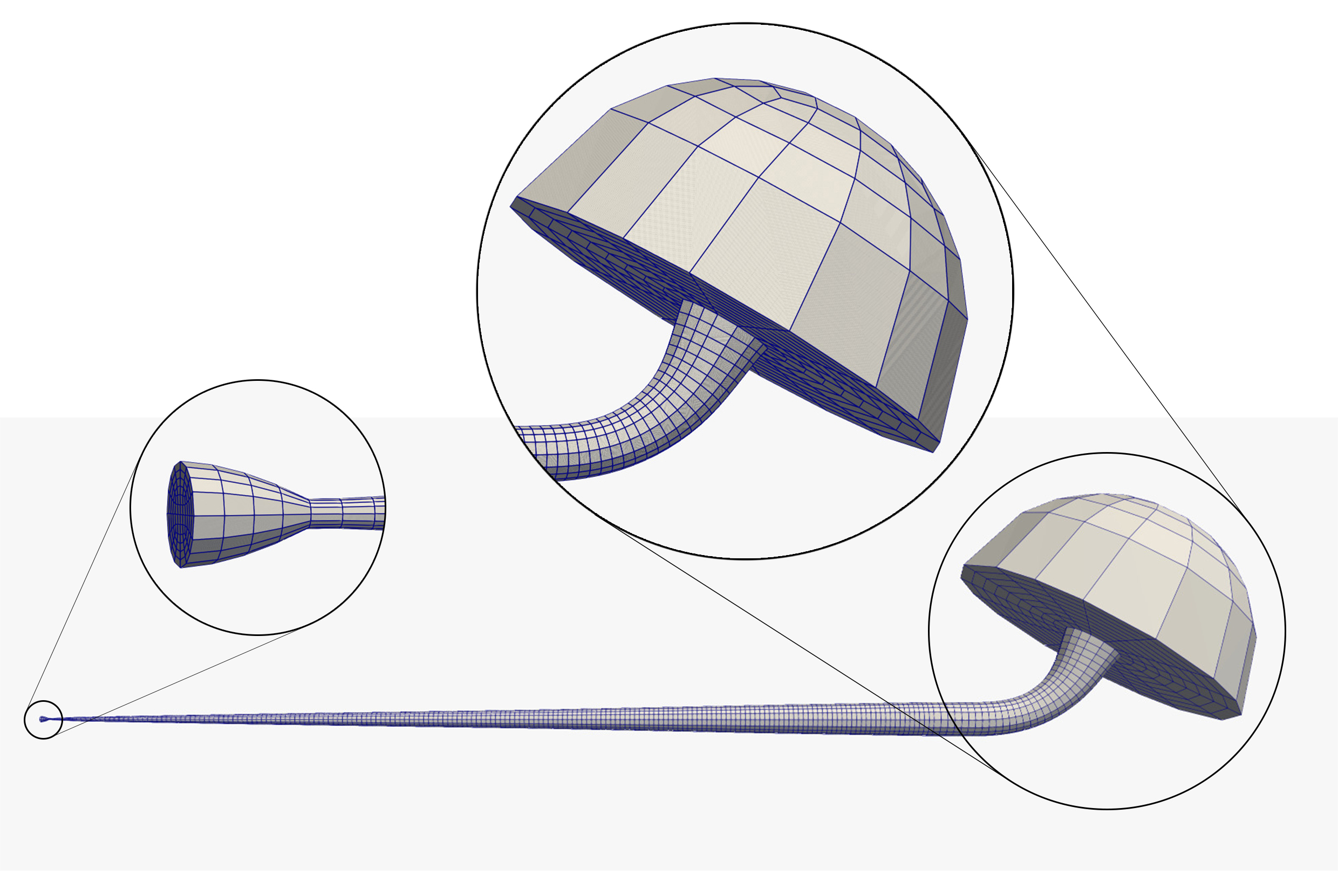}
        \caption{Volumetric mesh of the alphorn with a zoom on the mouthpiece and bell/hemispherical regions.}
        \label{fig:pre-processing:steps3}
    \end{subfigure}
    \caption{Overview of the reconstruction pipeline and final 3D mesh.}
    \label{fig:pre-processing}
\end{figure}

As reported in Fig.\ \ref{fig:alphorn_division}, the Swiss alphorn can be divided into 6 main chunks:
\begin{itemize}
    \item The mouthpiece, its reconstruction is based on the measurements of the piece;
    \item The four conical pieces. The radii and lengths are the parameters that vary in this work;
    \item The bell. Its reconstruction is based on the original mold (see Fig.\ \ref{fig:pre-processing:steps1}).
\end{itemize}
All the mentioned radii are given from direct measurements.
The centerline of the alphorn, except for the bell, is a straight line, where the length and the centerline position of each chunk is given by direct measurements. While for the bell, only the starting position (the one corresponding to $\mathrm{R_{5}}$) is available, conversely, the position related to $R_{bell}$ is not given, and it needs to be reconstructed.

The mouthpiece is composed of two different subparts: 
the first part (between $\mathrm{R_{inlet}}$ to $\mathrm{R_0}$) is non linear where the value of the radius for some points of the centerline are available from experimental measurements.
While the second one, between $\mathrm{R_0}$ and $\mathrm{R_1}$, is a conical piece. 

For the four big conical pieces the lengths ($\mathrm{L_1}, \mathrm{L_2}, \mathrm{L_3}, \mathrm{L_4}$) and the radii ($\mathrm{R_1}, \mathrm{R_2}, \mathrm{R_3}, \mathrm{R_4}$) are measured. Each piece described in Fig.\ \ref{fig:alphorn_division} is made separately, the alphorn then can be easily assembled. 

The bell requires a different treatment.
$R_5$, $R_{bell}$ and the starting point of the bell are available measurements. No information regarding the centerlines and the radii in between are given.
Due to its complex flaring profile and curvature, image processing techniques are applied to a photograph of the wooden mold (see Fig.\ \ref{fig:pre-processing:steps1}) used to build the bell. In particular, we first binarize the photograph, and then extract the contour (see Fig.\ \ref{fig:pre-processing:steps1}), then spline curves representing the centerline and the radius of the bell are computed.
Due to the limited resolution of the image, the extracted data have some noise. To replicate the smoothness of the instrument, we impose monotonicity to the radius data, while performing an isotonic regression of the centerline data. The results are again reported in Fig.\ \ref{fig:pre-processing:steps2}. 

Finally, the centerlines and radii of the six parts are separately (to avoid oscillation) interpolated with cubic splines, densely evaluated and concatenated to generate a unique parametrization of the alphorn. These data feed the mesh generator presented in \cite{MARCINNO2025118153} for generating the volumetric mesh reported in Fig.\ \ref{fig:pre-processing:steps3}. 
It is worth noting that the hemispherical region downstream of the bell must be constructed to correctly impose the boundary condition and avoid spurious reflections (see Sec.\ \ref{sec:mathematical_modeling}). The mesh generator handles the challenging task of constructing elements in the transition region between the bell exit and the hemispherical domain. This latter is constructed with a radius consistent with \cite{fehlmann1994} for our frequency range, set to $R = 0.447$ m.

\section{The Mathematical Framework}\label{sec:mathematical_modeling}
Sound is physically modeled as acoustic waves, which are small oscillations of pressure in a compressible acoustic medium, namely the air. 
We decide to model the sound using the stationary linear acoustic equation in the frequency domain, the Helmholtz equation.
Through this work, we fix the speed of the sound in the air $c_0 = 340 \,\,\mathrm{m\,s^{-1}}$ and the air density $\rho = 1.225 \,\, \mathrm{kg\, m^{-3}}$.

\subsection{The Helmholtz equation and the boundary conditions}
\begin{figure}
    \centering
    \includegraphics[width=\linewidth]{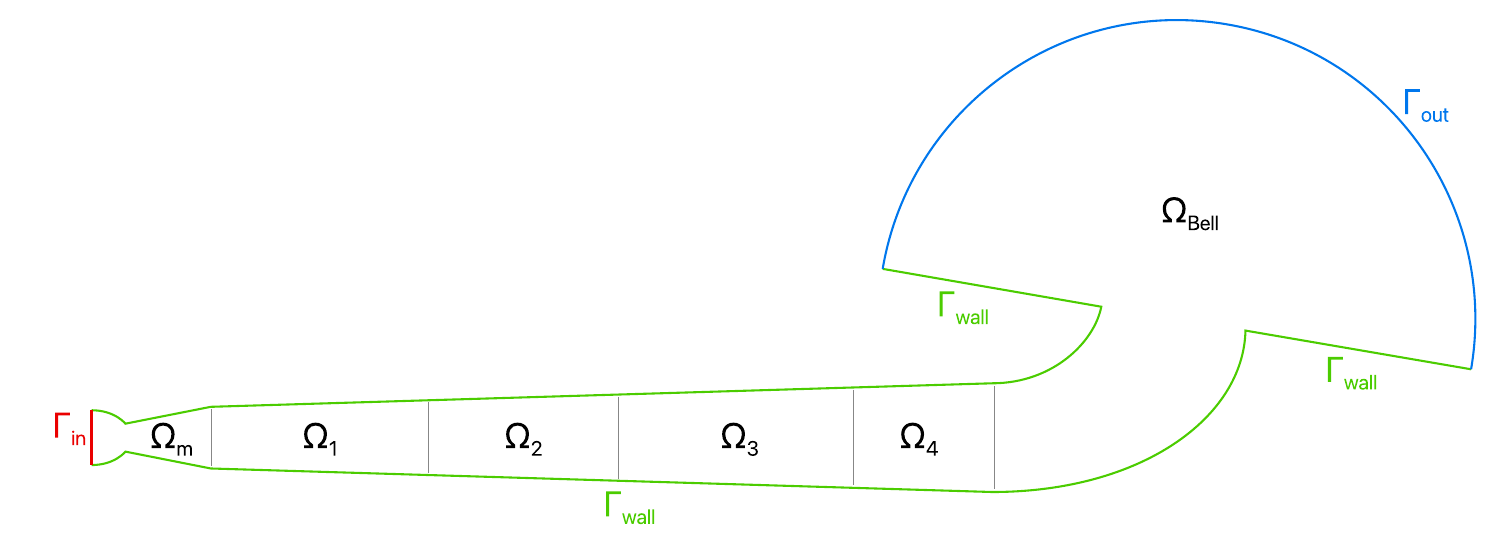}
    \caption{2D section of domains and boundaries. The colors are intended to help the comprehension of the boundary subdivisions: red for the inlet section, green for the wall boundary, and blue for the outlet boundary.}
    \label{fig:domain_alphorn}
\end{figure}

Referring to Fig.\ \ref{fig:domain_alphorn}, we define $\Omega$ as the total volume of the alphorn, where $\Omega = \Omega_m \cup \Omega_1  \cup \Omega_2 \cup
\Omega_3 \cup \Omega_4 \cup \Omega_{bell}$. $\Omega_{bell}$ represents both the bell and hemisphere domains.
$\Gamma_{in}$ is the inlet section of the mouthpiece, $\Gamma_{out}$ is the outlet of the hemispherical part and  $\Gamma_{wall}$ is the boundary of the conical and baffle pieces.
Let $p(\mathbf{x})$ denote the acoustic pressure field defined on $\Omega$, the interior of the instrument. 
The Helmholtz equation states \cite{ihlenburg1998finite}:
\begin{align*}
\Delta p + k^2 p &= 0 \quad \text{in} \quad \Omega \,,
\end{align*}
where $k = \nicefrac{2\pi f}{c_0}$ is the wave number and $f$ is the imposed frequency. The equation needs to be equipped with proper boundary conditions.

Firstly, we model the lateral walls of the instrument as "sound-hard" solid walls, therefore applying homogeneous Neumann boundary conditions:
\begin{equation*}
    \frac{\partial p}{\partial \mathbf{n}} = 0 \quad  \text{on} \quad \Gamma_{wall} \,,
    \label{eq:wall_bc}
\end{equation*}
where $\mathbf{n}$ is the outward normal vector to the surface.

The second boundary condition is related to wave generation at the inlet of the instrument. Although the interaction between the lips and the mouthpiece is still an open research field, a good first approximation, discussed in \cite{martin1942} and still used today in works like \cite{resch2019}, consists in considering the lip motion as sinusoidal, with a time harmonic velocity $U(\mathbf{x},t) = \text{Re}\{u_n e^{-i\omega t}\}$ where the constant $u_n$ is the imposed constant normal velocity at the mouthpiece. This results in an inhomogeneous Neumann boundary condition:
\begin{equation*}
\label{eq:inlet_bc}
    \frac{\partial p}{\partial \mathbf{n}} = i \rho \omega u_n \quad  \text{on} \quad \Gamma_{in} \,. \\
\end{equation*}

Lastly, we need to consider how the sound behaves at the outlet of the instrument. Since sound propagates outward from the bell into open space, we must artificially extend the domain beyond the instrument itself. As mentioned in Sec.\ \ref{sec:preprocessing_and_mesh_generation}, a common and effective approach is to enclose the bell extremity in a hemispherical domain, representing a portion of the surrounding free field. It is important to note that this hemisphere is not a physical boundary, but a truncation of the infinite domain. To avoid artificial reflections and accurately reproduce the effect of outgoing waves, an impedance boundary condition (IBC) is applied on the spherical surface of the hemisphere. 
As the name suggests, the IBC relates the pressure and the normal component of the velocity at the boundary through a prescribed surface impedance. In this case, we use the specific impedance for spherical wave propagation in air $Z = Z_{sph}$, which ensures that the boundary mimics the behavior of sound radiating freely into open space. Mathematically, this is expressed with \cite{fehlmann1994}:
\begin{equation*}
\label{eq:impedance}
    p = Z_{sph} u \quad \text{with} \quad Z_{sph}(r, k) = \rho c_0\frac{i k r}{1+ikr} \,,
\end{equation*}

depending on the radial distance from the source $r$ and the wave number $k$.
We obtain an inhomogeneous Robin boundary condition of the form:

    \begin{equation*}
    \label{eq:impedance_bc}
        \frac{\partial p}{\partial \mathbf{n}} = i \omega \frac{(1+ikr)}{c_0(i k r)}  p \, \quad  \text{on} \quad \Gamma_{out} \,.
    \end{equation*}

The final strong formulation along with the boundary conditions reads:\\

Find $p(\mathbf{x}) \in \mathcal{C}^2(\Omega;\mathbb{C})$ satisfying:
\begin{align}
\Delta p + k^2 p &= 0 \quad &&\text{in } \Omega \,, \label{eq:helmholtz}\\
\frac{\partial p}{\partial \mathbf{n}} &= i \rho \omega u_n \quad &&\text{on } \Gamma_{in} \,, \label{eq:inflow}\\
\frac{\partial p}{\partial \mathbf{n}} &= 0 \quad &&\text{on } \Gamma_{wall} \,, \label{eq:wall}\\
\frac{\partial p}{\partial \mathbf{n}} &= i \omega \frac{(1+ikr)}{c_0(i k r)}p \quad &&\text{on } \Gamma_{out} \,. \label{eq:outflow}
\end{align}

We proceed by deriving the weak formulation essential for the construction of the FEM. After having introduced the Sobolev space:

    \begin{equation*}
        H^1(\Omega; \mathbb{C}) = \left\{ v \in L^2(\Omega; \mathbb{C}) : \nabla v \in L^2(\Omega;\mathbb{C}) \right\} \,,
    \end{equation*}

we multiply Eq.\ \ref{eq:helmholtz} by a test function $v \in H^1(\Omega;\mathbb{C})$ and integrate over the domain $\Omega$:

    \begin{equation*}
        \int_{\Omega} (\Delta p + k^2 p) v = 0 \,.
    \end{equation*}

Following the standard procedure and integrating by parts (see \cite{quarteroni1994numerical}), we obtain the following weak form:

\begin{align}
    \left\{
    \begin{aligned}
    &\text{Find p} \in H^{1}(\Omega;\mathbb{C}) \,\, \text{such that:}\\
    &a(p,v) = F(v), \qquad \forall v \in H^{1}(\Omega;\mathbb{C}),
    \end{aligned}
    \right.
    \label{eq:weak_form}
\end{align}

where, the bilinear form $a(\cdot,\cdot)$ is defined as follows:

\begin{equation}
a(p,v) \coloneqq \int_{\Omega} \nabla p \cdot \nabla v
- k^{2} \int_{\Omega} p\,v
- \frac{1 + i k r}{r} \int_{\Gamma_{\mathrm{out}}} p\,v ,
\label{eq:bilinear_form}
\end{equation}

and the linear functional $F(\cdot)$ as:
\begin{equation}
F(v) \coloneqq i \rho \omega \int_{\Gamma_{\mathrm{in}}} u_n\, v .
\label{eq:linear_form}
\end{equation}

\subsection{The parametrization of the equation}
Solving Eq.\ \ref{eq:weak_form} using the FEM for a large number of geometries is unfeasible due to the high computational costs. With the idea of finding the best shape for the alphorn, we use the
reduced basis method (RBM) \cite{hesthaven2016certified}.
The RBM is a model order reduction technique for
parametrized partial differential equations. It relies on an offline-online decomposition: in the offline phase, high-fidelity solutions (called snapshots) are computed for selected parameter values and used to build a low-dimensional approximation space (the reduced basis). In the online phase, this reduced basis is used to rapidly compute approximate solutions for new parameter values.

The model introduced before is defined on a fixed domain $\Omega$. However, our goal is to simulate
different instrument shapes within the constrains of the manufacturer, which requires considering a family of domains $\Omega_\mu$
parametrized by a geometrical parameter $\mu$ in a parameter space
$\mathcal{P} \subset \mathbb{R}^d$.
To apply the RBM, we fix a reference parameter $\hat{\mu} \in \mathcal{P}$ and define a reference
domain $\hat{\Omega} := \Omega_{\hat{\mu}}$.
The problem is then reformulated using a bijective and sufficiently smooth mapping
$T_\mu : \hat{\Omega} \rightarrow \Omega_\mu$, which leads to the parametrized weak formulation.
Following Fig.\ \ref{fig:alphorn_division}, an alphorn is composed of a mouthpiece, four cones and a bell where the parametrization fixes the total length $L$ of the four conical pieces, as well as the geometry of the mouthpiece (both its position and radii, $R_{in}, R_{0}, R_{1}$) and of the bell and the hemisphere (both position and radii, $R_5$ and $R_{bell}$). 
Therefore, the geometry is parametrized using a six-dimensional parameter space:
\[
\mu = (l_1, l_2, l_3, R_2, R_3, R_4),
\]
where the length parameters are defined as normalized quantities
$l_i = L_i / L$.
We add constraints on the parameters to ensure physically realistic and constructible
shapes. The lengths $l_i$ must satisfy $l_i \geq l_{\min}$. The radii $R_i$ are increasing, i.e.\ $R_i \leq R_{i+1}$ for $i = 1, \dots, 4$. The values of $l_{\min}$, $R_1$, and $R_5$
are chosen based on experimental measurements reflecting practical manufacturing constraints.
The minimum cone length is set to
$L_{\min} = 220~\mathrm{mm}$. 
The minimum radius is set to $R_{\min} = 6.5~\mathrm{mm}$,
which corresponds to the smallest radius that can be reliably constructed, while the maximum radius is fixed to
$R_{\max} = 51.5~\mathrm{mm}$, so that the conical part connects consistently to the bell geometry.

For a given geometry $\mu \in \mathcal{P}$, following Fig.\ \ref{fig:domain_alphorn}, let us denote by
$\Omega_{\mathrm{m}}$, $\Omega_{\mu,\mathrm{1}}$,
$\Omega_{\mu,\mathrm{2}}$, $\Omega_{\mu,\mathrm{3}}$,
$\Omega_{\mu,\mathrm{4}}$, and $\Omega_{\mathrm{bell}}$
the different subvolumes of the domain. The parameters $\mu$ do not affect the mouthpiece and the bell. Setting
\begin{equation}
    b(p,v,\Omega) = \int_{\Omega} \nabla p \cdot \nabla v
- k^{2} \int_{\Omega} p\,v,
\label{eq:bilinear-to-parametrize}
\end{equation}

the bilinear form (see Eq.\ \ref{eq:bilinear_form}) of the weak formulation is rewritten for a generic parametrized domain $\Omega_\mu$ as follows:
\begin{equation}
a(p,v) =
b(p,v,\Omega_{\mathrm{m}})
+ b(p,v,\Omega_{\mathrm{bell}})
+ \sum_{c=1}^{4} b(p,v,\Omega_{\mu,\mathrm{c}})
- \frac{1 + i k r}{r} \int_{\Gamma_{\mathrm{out}}} p\,v ,
\label{eq:parametrized_weak_form}
\end{equation}

where the subscript $c$ denotes the index of the conical piece.
The right-hand side term, $F(v)$, is not affected by the parametrization.


\subsubsection{Affine decomposition of the parametric form}\label{sec:affine_decomposition}
To remove the dependence on $\mu$ from the integration domains appearing in the terms
$\sum_{c=1}^{4} b(p,v,\Omega_{\mu,\mathrm{i}})$ in Eq.\ \ref{eq:parametrized_weak_form}, we introduce a change of variables on a generic conical segment.
We consider an abstract cone, denoted by the subscript $c$, whose geometry is described by the
four parameters:
\[
\mu_c \coloneqq (x_1,x_2,r_1,r_2) \, ,
\]

where $x_1$ and $x_2$ are the positions of the initial and ending sections, respectively, and $r_1$ and $r_2$ are the radii of the initial and ending sections, respectively.
The generic cone parameterization is defined as:
\[
\Omega_{\mu_c} \coloneqq \left\{(x,y,z)\in\mathbb{R}^3 \,\middle|\, 
x\in[x_1,x_2],\;
\sqrt{y^2+z^2}\le r_1 + \frac{r_2-r_1}{x_2-x_1}(x-x_1)
\right\}.
\]

We fix a reference set of parameters and the associated reference cone:
\[
\hat{\mu}_c \coloneqq (\hat{x}_1,\hat{x}_2,\hat{r}_1,\hat{r}_2)\,, \quad 
\hat{\Omega}_c \coloneqq \Omega_{\hat{\mu}_c}.
\]
Our goal is to rewrite the integrals $\int_{\Omega_{\mu_c}} \nabla p\cdot\nabla v$ and
$\int_{\Omega_{\mu_c}} pv$ as integrals over $\hat{\Omega}_c$, moving outside the integral the dependence from the geometrical parameters.
To this end, we construct a bijective and sufficiently smooth mapping $ T_{\mu_c} : \hat{\Omega}_c \to \Omega_{\mu_c}$ obtained as a combination of an axial rescaling and a radial
stretching,  defined as follows:

\begin{equation}
T_{\mu_c}(x,y,z) \coloneqq
\bigl(x_1 + \alpha(\mu_c)(x-\hat{x}_1),\;
\lambda(x;\mu_c)\,y,\;
\lambda(x;\mu_c)\,z\bigr),
\label{eq:mapping}
\end{equation}

where:
\[
\alpha(\mu_c) \coloneqq \frac{x_2-x_1}{\hat{x}_2-\hat{x}_1},
\qquad
\lambda(x;\mu_c) \coloneqq \frac{N(x;\mu_c)}{D(x)},
\]

and consequently:

\[
N(x;\mu_c) \coloneqq r_1 + \beta(\mu_c)(x-\hat{x}_1),
\qquad
D(x) \coloneqq \hat{r}_1 + \hat{\gamma}(x-\hat{x}_1),
\qquad
\beta(\mu_c) \coloneqq \frac{r_2-r_1}{\hat{x}_2-\hat{x}_1},
\qquad
\hat{\gamma} \coloneqq \frac{\hat{r}_2-\hat{r}_1}{\hat{x}_2-\hat{x}_1}.
\]

For any function $f$ defined on the parametric cone
$\Omega_{\mu_c}$, we define its pullback to the reference cone $\hat{\Omega}_c $ with
\[
\hat{f} \coloneqq f \circ T_{\mu_c}.
\]
Applying the chain rule yields
\[
\nabla \hat{f}(\xi)
= D T_{\mu_c}(\xi)^{T}\,\nabla f\bigl(T_{\mu_c}(\xi)\bigr),
\]
and therefore
\[
\nabla f\bigl(T_{\mu_c}(\xi)\bigr)
= \bigl(D T_{\mu_c}(\xi)^{T}\bigr)^{-1}\nabla \hat{f}(\xi).
\]

The computation of the Jacobian matrix $(D T_{\mu_c}(\xi))$ is detailed in the Appendix \ref{app:appendix_1}.
Using the change of variables $\xi = T_{\mu_c}(\hat{\xi})$, the bilinear form on the cone can be
rewritten as
\[
\int_{\Omega_{\mu_c}} \nabla p \cdot \nabla v
= \int_{\hat{\Omega}_c} (\nabla \hat{p})^{T} G_{\mu_c}\,\nabla \hat{v},
\]
and
\[
\int_{\Omega_{\mu_c}} pv
= \int_{\hat{\Omega}_c} \hat{p}\,\hat{v}\,\det\bigl(DT_{\mu_c}\bigr).
\]

The parameter dependence now appears only through $G_{\mu_c}$ and
$\det\bigl(DT_{\mu_c}\bigr)$ (see Appendix \ref{app:appendix_2} and \ref{app:appendix_3} for their definitions), representing the metric tensor and the determinant of the Jacobian matrix, respectively. To achieve an efficient implementation of the RBM, these quantities are expressed as affine decompositions with respect to $\mu_c$, respectively:
\[
G_{\mu_c}(\xi) = \sum_{i=1}^{7} \omega_i(\mu_c)\, G_i(\xi),
\qquad
\det\bigl(DT_{\mu_c}(\xi)\bigr)
= \sum_{j=1}^{3} \theta_j(\mu_c)\, q_j(\xi).
\]

The explicit construction of the affine decomposition of the determinant of the Jacobian and the metric tensor is given in the Appendix \ref{app:appendix_2} and \ref{app:appendix_3}.
As a result, the bilinear forms on the cone take the form:
\begin{equation}
\int_{\Omega_{\mu_c}} \nabla p \cdot \nabla v
= \sum_{i=1}^{7} \omega_i(\mu_c)
\int_{\hat{\Omega}_c} (\nabla \hat{p})^{T} G_i\, \nabla \hat{v}, \quad
\int_{\Omega_{\mu_c}} pv
= \sum_{j=1}^{3} \theta_j(\mu_c)
\int_{\hat{\Omega}_c} q_j(\xi)\, \hat{p}\,\hat{v}.
\label{eq:param2}
\end{equation}

Combining the previous findings, we obtain an affine parametrization of the weak formulation defined on the entire reference domain.
The bilinear form can therefore be written as:
\begin{align*}
a(p,v) ={}&
b(p,v,\Omega_{\mathrm{m}})
+ b(p,v,\Omega_{\mathrm{bell}}) + \\
& \sum_{c=1}^{4} \left(
    \sum_{i=1}^{7} \omega_i(\mu_c)
    \int_{\hat{\Omega}_c} (\nabla \hat{p})^{T} G_i\, \nabla \hat{v}
    + \sum_{j=1}^{3} \theta_j(\mu_c)
    \int_{\hat{\Omega}_c} q_j(\xi)\, \hat{p}\,\hat{v}
\right) - \\
& \frac{1 + i k r}{r} \int_{\Gamma_{\mathrm{out}}} p\,v .
\end{align*}

To numerically solve the problem with the FEM, we discretize the problem using the Galerkin method \cite{quarteroni1994numerical}, and then the RBM is applied thanks to the affine decomposition of the matrices.
In particular, this decomposition yields the standard \emph{offline--online} structure, where all the parameter-independent quantities are assembled once in an \emph{offline} stage, while the
\emph{online} stage only requires the evaluation of the parameter-dependent coefficients $\omega_i$ and $\theta_j$ and the solution
of a reduced linear system \cite{hesthaven2016certified}.

\subsection{The final pipeline}
\begin{figure}[h!]
    \centering
        \includegraphics[width=\linewidth]{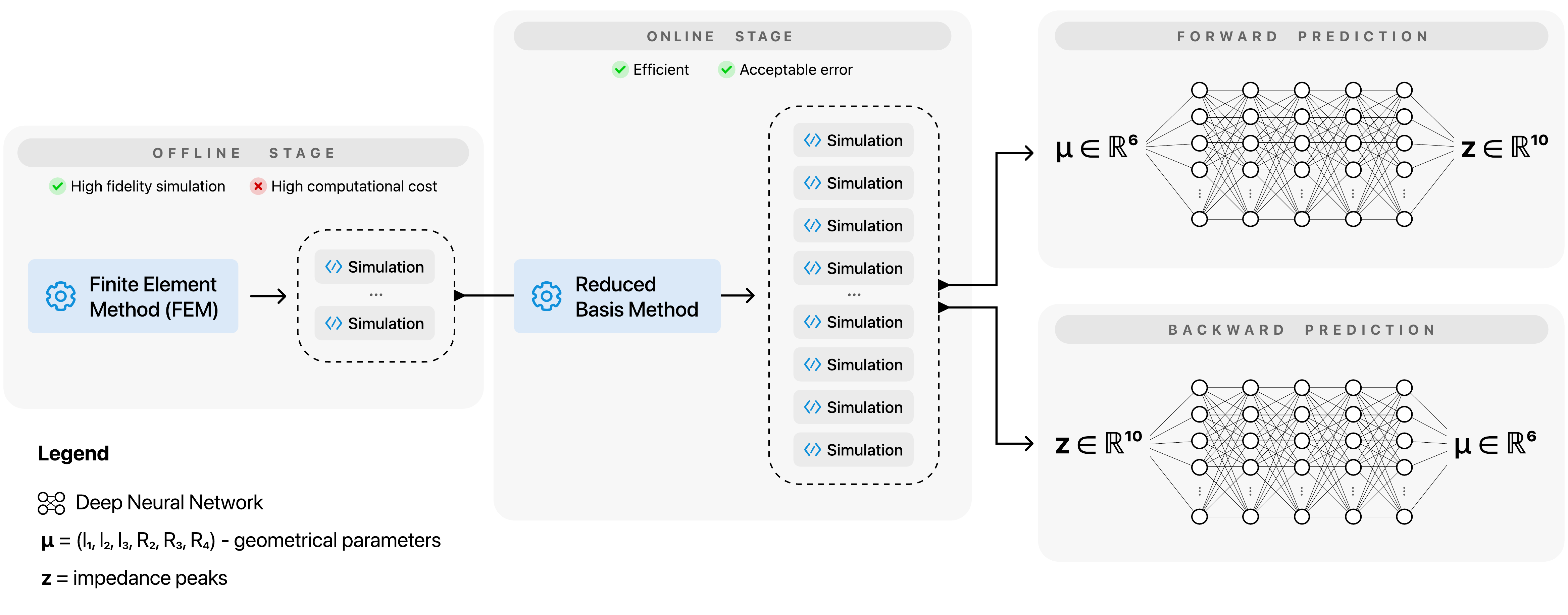}
    \caption{Final pipeline of the optimization process.}
    \label{fig:pipeline}
\end{figure}

Referring to Fig.~\ref{fig:pipeline}, we present the pipeline for the shape optimization of the
alphorn. In particular, the procedure consists of the following steps:
\begin{itemize}
    \item Generation of FEM simulations using the parametrized equations (this avoids the
    generation of a new mesh);
    \item Generation of thousands of simulations through the reduced basis framework;
    \item Use of a deep neural networks to identify the optimal alphorn shape.
\end{itemize}

\section{Numerical Results}\label{sec:results}
Sec.~\ref{sec:experimental_alphorns} presents a comparison of three experimentally measured alphorns and defines the reference instrument used in the remainder of the study.
Sec.~\ref{sec:basis_construction} details the generation of the reduced basis and, consequently, the dataset of simulations. 
Sec.~\ref{sec:optimization_results} discusses the architectures of the ML models and their main numerical results.
Appendix~\ref{apx:verification_on_test} provides a validation of the FEM implementation.

To ensure an efficient RB framework, we adopt a mesh with an average cell volume of approximately $0.010 \text{cm}^3$, also accounting for the inherently larger cells in the hemisphere (see Fig.~\ref{fig:pre-processing}). A convergence study performed on the alphorn mesh size shows that the resulting impedance spectrum remains essentially unchanged, with only minor deviations compared to results obtained using finer meshes.

\subsection{Selection of the reference instrument}
\label{sec:experimental_alphorns}
Three real alphorns were experimentally measured to evaluate the variability of their acoustic response and to establish a reliable reference configuration. All the three instruments share the same mouthpiece and bell geometries. The different geometrical parameters, the ones defining the four conical pieces, are reported in Appendix \ref{apx:real-alphorns-measurements}.

For each instrument, seven playable notes were recorded in order to produce predefined target pitches ($f_{\mathrm{note}}$). Several experimental attempts were performed for each note, and the frequency corresponding to the best performance, defined as the closest to the target pitch, was retained as the measured frequency ($f_{\mathrm{meas}}$).

\begin{figure}[t]
    \centering
    \includegraphics[width=\linewidth]{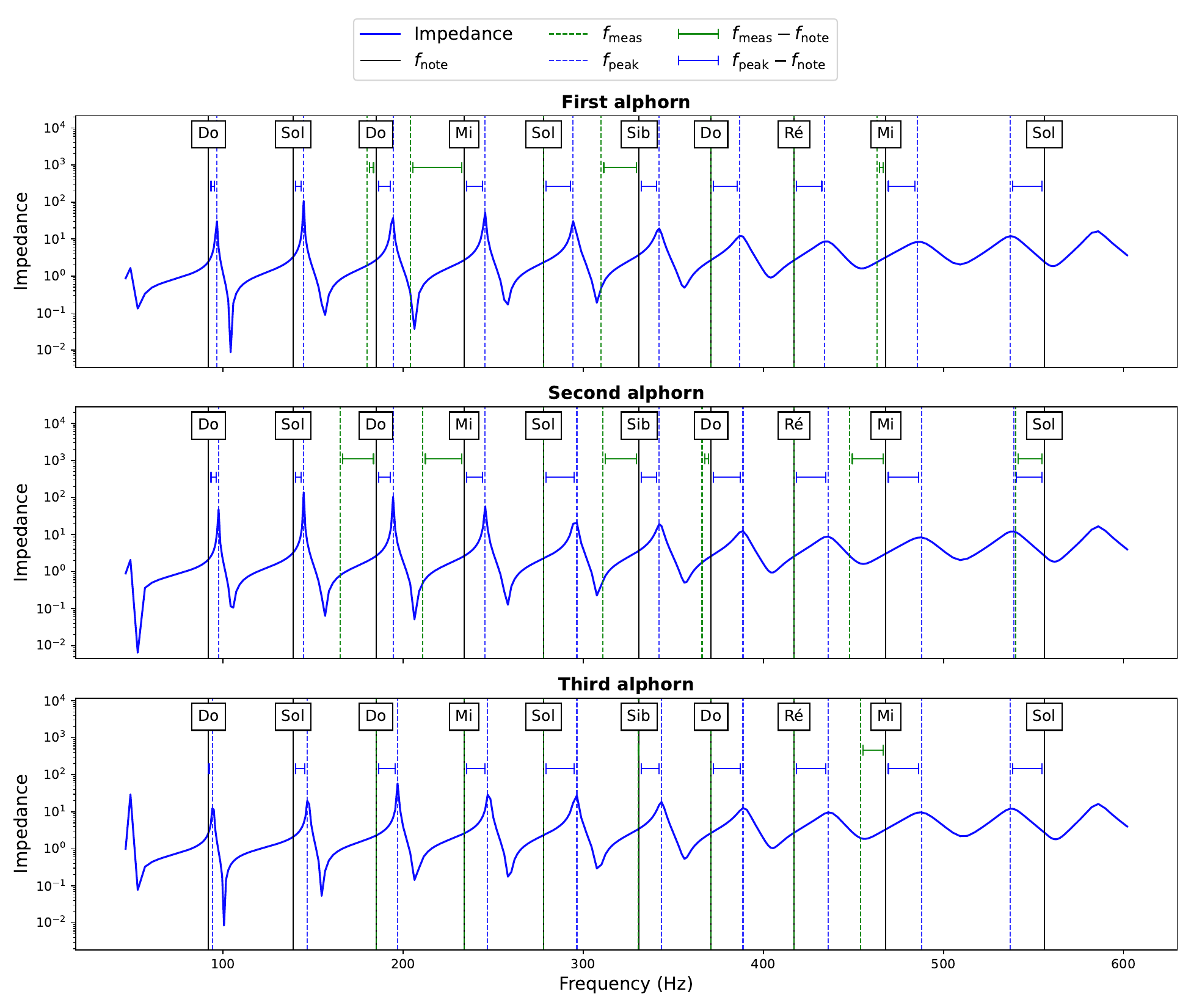}
    \caption{Simulated impedance of the three alphorns as a function of frequency. Black vertical lines indicate frequencies of the alphorn scale $f_{\mathrm{note}}$, green dashed vertical lines are the measured frequencies $f_{\mathrm{meas}}$, and blue dashed lines are the resonance frequencies $f_{\mathrm{peak}}$ extracted from FEM impedance spectra. Horizontal markers show the corresponding deviations with respect to $f_{\mathrm{note}}$.}
    \label{fig:real_alphorns}
\end{figure}

The impedance spectra of the three alphorns were also simulated using FEM. From these spectra, the corresponding numerical resonance frequencies ($f_{\mathrm{peak}}$) were extracted.
The results are reported in Tab. \ref{tab:real_alphorns} and Fig.~\ref{fig:real_alphorns}.

\renewcommand{\arraystretch}{1.3}

\begin{table}[H]
\centering

\begin{tabular}{lccccccc}

\multicolumn{8}{c}{\textbf{First alphorn}} \\
\hline
Notes & Do$_1$ (C) & Mi$_1$ (E) & Sol (G) & Sib (B$\flat$) & Do$_2$ (C) & Re (D) & Mi$_2$ (E) \\
\hline
Target frequency $f_{\mathrm{note}}$ (Hz) & 185.00 & 234.00 & 278.00 & 331.00 & 371.00 & 417.00 & 468.00 \\
Measured frequency $f_{\mathrm{meas}}$ (Hz) & 180.00 & 204.00 & 278.00 & 310.00 & 371.00 & 417.00 & 463.00 \\
Measured deviation $\Delta f_{\mathrm{meas}}$ (cents) & -47.43 & -237.53 & +0.00 & -113.48 & +0.00 & +0.00 & -18.60 \\
FEM frequency $f_{\mathrm{peak}}$ (Hz) & 194.47 & 245.60 & 294.42 & 342.20 & 386.83 & 433.89 & 485.48 \\
Deviation $\Delta f_{\mathrm{peak}}$ (cents) & +86.42 & +83.75 & +99.35 & +57.62 & +72.33 & +68.72 & +63.49 \\
\hline \\

\multicolumn{8}{c}{\textbf{Second alphorn}} \\
\hline
Notes & Do$_1$ & Mi$_1$ & Sol & Sib & Do$_2$ & Re & Mi$_2$ \\
\hline
Target frequency $f_{\mathrm{note}}$ (Hz) & 185.00 & 234.00 & 278.00 & 331.00 & 371.00 & 417.00 & 468.00 \\
Measured frequency $f_{\mathrm{meas}}$ (Hz) & 165.00 & 211.00 & 278.00 & 311.00 & 366.00 & 417.00 & 448.00 \\
Measured deviation $\Delta f_{\mathrm{meas}}$ (cents) & -198.07 & -179.12 & +0.00 & -107.90 & -23.49 & +0.00 & -75.61 \\
FEM frequency $f_{\mathrm{peak}}$ (Hz) & 194.47 & 245.60 & 296.51 & 342.20 & 388.71 & 436.02 & 487.77 \\
Deviation $\Delta f_{\mathrm{peak}}$ (cents) & +86.42 & +83.75 & +111.57 & +57.62 & +80.74 & +77.22 & +71.62 \\
\hline\\

\multicolumn{8}{c}{\textbf{Third alphorn}} \\
\hline
Notes & Do$_1$ & Mi$_1$ & Sol & Sib & Do$_2$ & Re & Mi$_2$ \\
\hline
Target frequency $f_{\mathrm{note}}$ (Hz) & 185.00 & 234.00 & 278.00 & 331.00 & 371.00 & 417.00 & 468.00 \\
Measured frequency $f_{\mathrm{meas}}$ (Hz) & 185.00 & 234.00 & 278.00 & 330.50 & 371.00 & 417.00 & 454.00 \\
Measured deviation $\Delta f_{\mathrm{meas}}$ (cents) & +0.00 & +0.00 & +0.00 & -2.62 & +0.00 & +0.00 & -52.58 \\
FEM frequency $f_{\mathrm{peak}}$ (Hz) & 196.97 & 246.98 & 296.51 & 343.50 & 388.71 & 436.02 & 487.77 \\
Deviation $\Delta f_{\mathrm{peak}}$ (cents) & +108.51 & +93.44 & +111.57 & +64.17 & +80.74 & +77.22 & +71.62 \\
\hline \\
\end{tabular}
\caption{Measured and FEM-derived resonance frequencies compared to the alphorn scale frequencies for the three alphorns. The resonance frequencies $f_{\mathrm{peak}}$ are extracted from FEM-computed impedance spectra. The deviations $\Delta f_{\mathrm{meas}}$ and $\Delta f_{\mathrm{peak}}$ are expressed in cents relative to $f_{\mathrm{note}}$.}

\label{tab:real_alphorns}

\end{table}

Among the three instruments, the third alphorn exhibits the best agreement between the measured frequencies and the frequencies of the alphorn scale. It is therefore taken as the reference domain in the remainder of the study, defining the total length $L$, as well as the geometry of the mouthpiece and the bell.

\subsection{Construction of the reduced basis \& dataset construction}\label{sec:basis_construction}
The reduced basis is constructed with a proper orthogonal decomposition (POD) strategy. A large dataset of high-fidelity solutions is generated by varying both the geometric parameter $\mu$ and the frequency $f$. 
In total, 1225 geometric configurations and 200 frequency values are considered, resulting in 245000 snapshots.
Since a direct singular value decomposition (SVD) of this dataset was computationally prohibitive, we employ a randomized SVD algorithm \cite{doi:10.1137/090771806}. This approach yielded satisfactory results at low frequencies, with a relative error of order $10^{-3}$ for frequencies below 200~Hz, while for higher frequencies the results were of order greater than $10^{-3}$, even with an increased number of POD modes. This indicates that a single
global reduced basis is not sufficient to accurately capture the solution behavior over the full
frequency range.
To address this limitation, frequency-localized reduced bases are introduced, guided by the
target notes one aims to produce on the alphorn, 92~Hz (Do), 139~Hz (Sol), 185~Hz (Do), 234~Hz (Mi), 278~Hz
(Sol), 331~Hz (Sib), 371~Hz (Do), 417~Hz (Ré), 468~Hz (Mi), and 556~Hz (Sol).
For each frequency zone, a dedicated reduced basis of dimension $N = 50$ is constructed
using snapshots restricted to the corresponding frequency interval. To ensure a smooth
transition between neighboring zones, the training intervals were chosen with deliberate
overlaps preventing a loss of accuracy near the boundaries of each zone.

The dataset of reduced numerical simulations is obtained as follows. 
For a given pair $(f,\mu)$, the frequency zone whose
center is closest to the target frequency $f$ is first identified. The corresponding reduced basis
is then selected, and the reduced system is assembled and solved within this space.
The frequency-localized reduced models are tested on random parameter values $\mu$ not
included in the snapshot set. The relative error between the reduced-order solution and the finite element solution is typically of order $10^{-3}$ to $10^{-2}$. 
Overall, this level of accuracy is considered satisfactory for the intended acoustic analysis,
given the strong reduction in computational cost.\footnote{This accuracy is achieved
with reduced spaces of dimension $N = 50$. In our experiments, the time required to solve the
reduced system was approximately 45--55~ms, compared to 2.4--2.6~s for the full finite
element solve, which corresponds to a measured speedup of about 50$\times$ (typically
50--56$\times$) for the same pairs $(f,\mu)$.}
Compared with the global reduced basis approach, which required reduced spaces of
dimension up to $N = 500$ without achieving uniform accuracy over the full frequency range,
the frequency-localized strategy offers a clear improvement in computational efficiency and precision.

Using the localized reduced bases, a dataset is generated by sampling $10^5$ different alphorn geometries across the admissible parameter space. For each geometry, the input impedance is evaluated over a dense frequency grid.
The frequency range of interest is partitioned into disjoint intervals, each associated with one target musical note. The boundaries of these intervals are defined as the midpoints
between consecutive target note frequencies. This construction yields a partition of the
spectrum with exactly one interval per note.
Within each interval, all peaks of the impedance magnitude are identified, and the one closest to the corresponding target note frequency is stored.
This selected peak is used as the resonance frequency associated with the note.
The resulting dataset therefore consists of tuples:
\[
(\boldsymbol{\mu}, f_1, \ldots, f_{10}),
\]
where $\boldsymbol{\mu}$ denotes the geometric parameters and $f_k$ are the extracted resonance frequencies, one for each target note.

\subsection{Results of the optimization}\label{sec:optimization_results}
The dataset of reduced basis simulations is used to train two machine learning models.
The first model addresses a forward problem: given a geometry $\boldsymbol{\mu}$, it predicts the associated
set of resonance frequencies $(f_1, \ldots, f_{10})$.
The second model addresses the inverse problem: given a target set of resonance frequencies, it predicts a geometry $\boldsymbol{\mu}$ expected to produce them. 

\subsubsection{Forward model}
The forward problem is treated as a supervised regression task with unconstrained outputs.
Input parameters and target frequencies are standardized, and the frequencies are represented
on a logarithmic scale.
The model is implemented as a multilayer perceptron (see Fig.\ \ref{fig:pipeline}) and trained using the Adam optimizer \cite{lecun2015deep}
with mini-batches of size $256$ and an $\ell_2$ regularization weight of $10^{-5}$.
Hyperparameters are selected using $4$-fold cross-validation (see Appendix~\ref{apx:hyper-forward}).
The final network architecture consists of five hidden layers with $256$ neurons each and ReLU activation functions. A learning rate of $10^{-3}$ is retained. With this configuration, the model achieves a mean squared error of $1.89 \times 10^{-4}$ on the
logarithm of the resonance frequencies.

\begin{figure}
    \centering
    \includegraphics[width=0.9\linewidth]{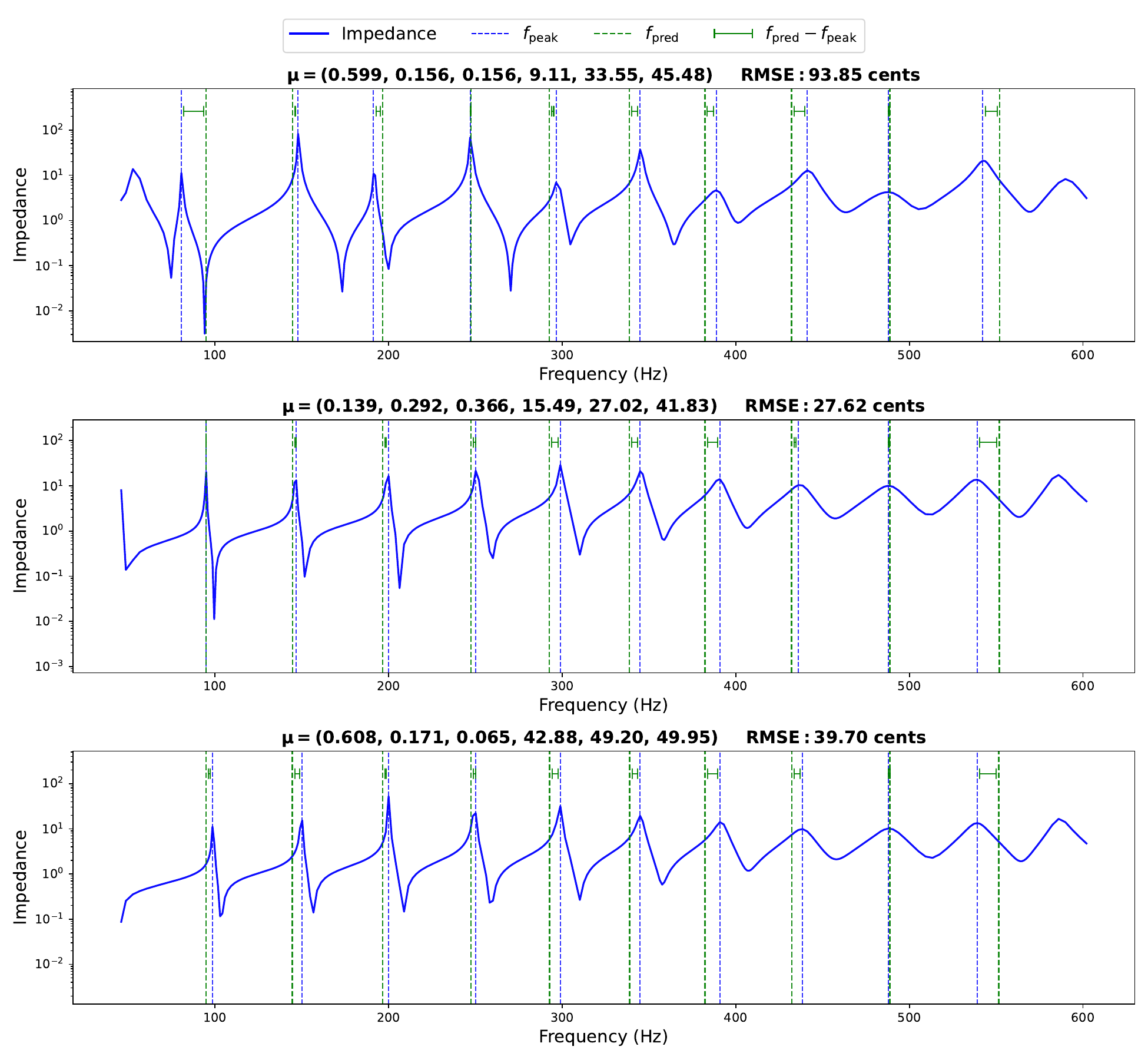}
    \caption{Impedance spectra for three randomly sampled geometries. Detected impedance peaks are indicated by dashed blue lines, while the resonance frequencies predicted by the ML forward model are shown in green. The discrepancy between the predicted resonance frequencies and those extracted from the impedance spectra is quantified using the root mean squared error (RMSE) of the deviations expressed in cents, after matching each predicted frequency to the closest impedance peak.}
    \label{fig:impedence_peak_result}
\end{figure}

Fig.~\ref{fig:impedence_peak_result} illustrates the performance of the forward model on three randomly selected geometries. 
The ML predicted resonance frequencies are compared with the peaks extracted from the corresponding high-fidelity FE impedance spectra.
They are generally close to the peaks of the impedance spectra, but are not perfectly aligned for all modes.
The error is quantified using the root mean squared error (RMSE) of the deviations expressed in cents, with values ranging from approximately 28 to 94 cents for the presented examples. 
These results suggest that the forward model can be used to approximate the resonance frequencies for large-scale exploration of the parameter space.

This approach is then used to perform a large-scale search over the admissible parameter space. This latter is built as a uniform grid of about $6\,500\,000$ possible geometrical parameters (the three lengths and radii of the cones), where a valid configuration can be extracted. 
For each configuration, the model predicts the associated resonance frequencies, which are compared to the target frequencies of the alphorn scale. 
Moreover, for each configuration, the mean squared error of the deviations in cents is computed and the configuration yielding the lowest error is retained.
The resulting geometric parameter vector, rounded to two decimal places, is
\[
\boldsymbol{\mu_{\text{best-forward}}} = (0.19,\; 0.09,\; 0.72,\; 6.50,\; 6.50,\; 29.90).
\]

\subsubsection{Backward model}\label{sec:backward-model}
The inverse mapping is inherently non-unique, as different geometries may lead to similar
acoustic responses. The objective is therefore to recover a physically admissible geometry
that matches a given set of target frequencies.
Geometric constraints are enforced directly at the architectural level, in fact,
the main challenge in designing the backward model is to ensure that the predicted geometric parameters satisfy the physical constraints. Specifically, the radii of the conical starting and ending sections have to be strictly increasing,
\begin{equation}
R_i < R_{i+1}, \quad \forall i \in \{1,2,3,4\},
\end{equation}
and the lengths must sum to one,
\begin{equation}
l_1 + l_2 + l_3 + l_4 = 1.
\end{equation}

We considered two approaches to enforce these constraints. The first approach consists of adding penalty terms to the loss function during training. However, this does not guarantee that the constraints are strictly satisfied; it merely encourages the model to satisfy them. As a result, the model may still produce invalid geometries, and the biased loss function may lead to suboptimal performance, so we discarded it.
The second (adopted) approach consists of designing the model architecture such that it can only output valid geometric parameters, thereby keeping the loss function unbiased. To enforce the length constraint, the model outputs four values corresponding to the segment lengths $(l_1, l_2, l_3, l_4)$. A softmax \cite{lecun2015deep} activation is then applied to these values, ensuring that they are all positive and sum to one. The fourth length $l_4$ is subsequently discarded.
Enforcing the radii constraints is more challenging. The boundary radii $R_1$ and $R_5$ are fixed, and strict monotonicity must be ensured for the intermediate radii. Rather than directly predicting the absolute values of $R_2$, $R_3$, and $R_4$, the model predicts the relative proportions of the intervals between consecutive radii.
Specifically, the model outputs the following four values corresponding to the intervals
\[
R_2 - R_1,\quad R_3 - R_2,\quad R_4 - R_3,\quad R_5 - R_4.
\]
A softmax activation is applied to these outputs, producing four positive values that sum to one. These values represent the relative size of each interval with respect to the total available radial growth $R_5 - R_1$.
By multiplying these values by $(R_5 - R_1)$, the physical size of each gap is obtained. The absolute values of the radii are then recovered by computing the cumulative sum of these gaps and adding the starting radius $R_1$. Since the gaps are strictly positive due to the softmax activation, the resulting radii are guaranteed to be strictly increasing and bounded between $R_1$ and $R_5$.
The hidden dimension is selected using a $4$-fold cross-validation procedure (see Appendix~\ref{apx:hyper-backward}), resulting in $1024$ neurons per layer. Each hidden layer uses a ReLU activation function. To reduce
overfitting, dropout is applied after each hidden layer with a rate of $0.1$. Training is
performed using the Adam optimizer \cite{lecun2015deep}, as in the forward model.

The backward model is evaluated by prescribing a set of target resonance frequencies
corresponding to the official notes of the alphorn scale. From these targets, the model
predicts a geometry expected to reproduce the desired acoustic response.
The resulting geometric parameter vector, rounded to two decimal places, is
\[
\boldsymbol{\mu_{\text{best-backward}}} = (0.28,\; 0.31,\; 0.26,\; 6.55,\; 11.85,\; 32.44).
\]

\subsubsection{Result discussion}
\begin{figure}[t]
    \centering
    \includegraphics[width=\linewidth]{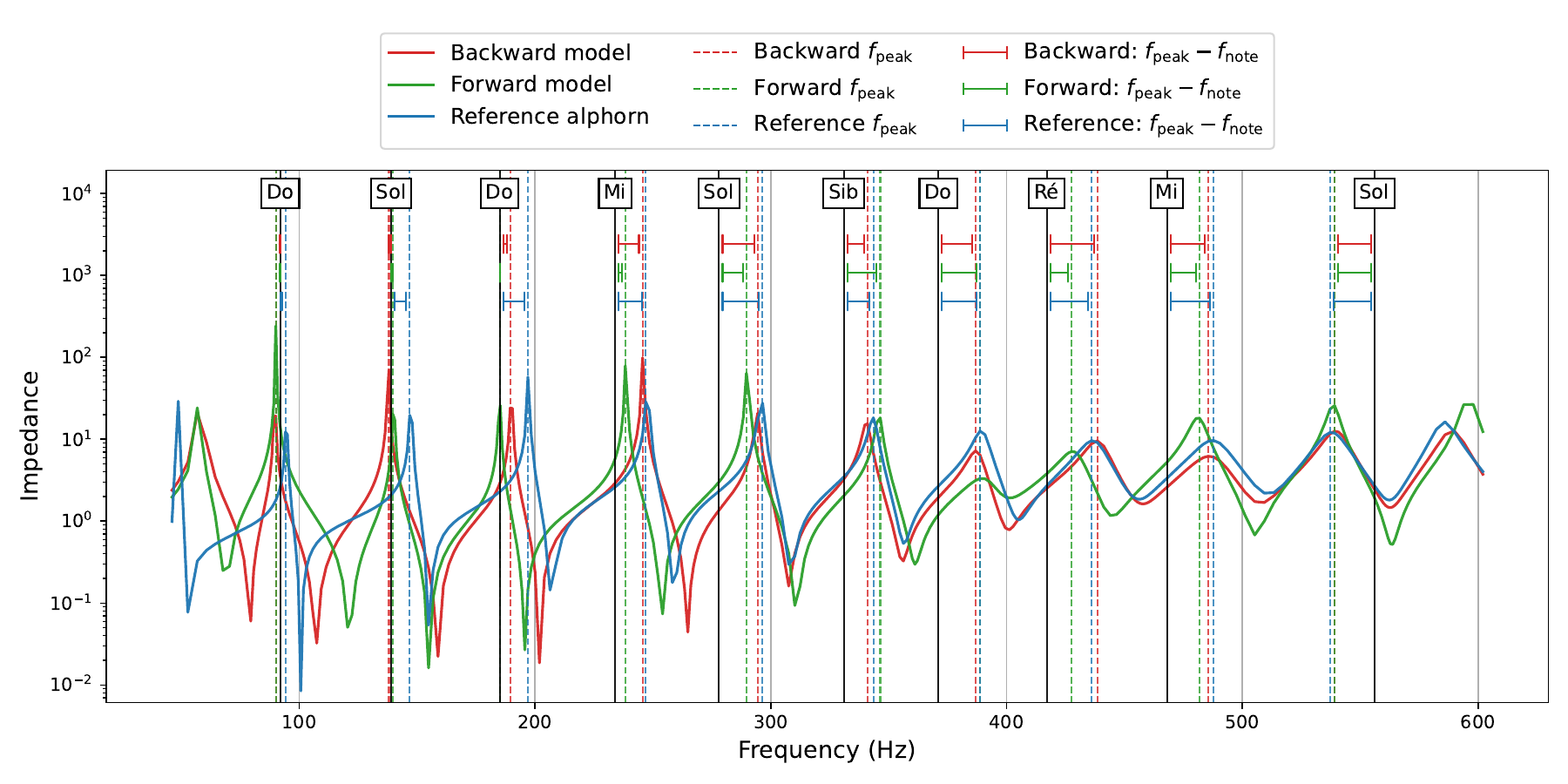}
    \caption{
    Comparison between the reference alphorn geometry and the geometries predicted by the forward ($\boldsymbol{\mu_{\text{best-forward}}}$) and backward ($\boldsymbol{\mu_{\text{best-backward}}}$) approaches. 
    For each case, the extracted FEM resonance frequencies $f_{\mathrm{peak}}$ are compared to the frequencies of the alphorn scale $f_{\mathrm{note}}$.}
    \label{fig:comparaison-back-forw-real}
\end{figure}

Fig.\ \ref{fig:comparaison-back-forw-real} and Tab.\ \ref{tab:final_comparison_models} present a direct visual comparison of the impedance spectra and resonance frequencies associated with the reference alphorn geometry (alphorn 3) and the best geometries predicted by the forward and backward ML models.
In particular, Tab.\ \ref{tab:final_comparison_models} reports the frequencies of the alphorn scale together with the resonance frequencies identified from the impedance spectra (simulated with FEM) of the reference and predicted geometries.
The deviations, expressed in cents, are also given for each note. The overall discrepancy is summarized using the root mean squared error (RMSE).

\begin{table}
\centering
\renewcommand{\arraystretch}{1.3}

\begin{tabular}{lcccccccccc}


\multicolumn{11}{c}{\textbf{Backward model}} \\
\hline
Notes & Do & Sol & Do & Mi & Sol & Sib & Do & Ré & Mi & Sol \\
\hline
$f_{\mathrm{note}}$ (Hz) & 92.00 & 139.00 & 185.00 & 234.00 & 278.00 & 331.00 & 371.00 & 417.00 & 468.00 & 556.00 \\
$f_{\mathrm{peak}}$ (Hz) & 90.01 & 138.03 & 189.47 & 245.60 & 294.42 & 341.00 & 386.83 & 438.49 & 485.48 & 539.11 \\
$\Delta f$ (cents) & -37.82 & -12.09 & +41.31 & +83.75 & +99.35 & +51.50 & +72.33 & +87.01 & +63.49 & -53.41 \\
\multicolumn{11}{r}{\textit{RMSE: 65.20 cents}} \\
\hline \\
\multicolumn{11}{c}{\textbf{Forward model}} \\
\hline
Notes & Do & Sol & Do & Mi & Sol & Sib & Do & Ré & Mi & Sol \\
\hline
$f_{\mathrm{note}}$ (Hz) & 92.00 & 139.00 & 185.00 & 234.00 & 278.00 & 331.00 & 371.00 & 417.00 & 468.00 & 556.00 \\
$f_{\mathrm{peak}}$ (Hz) & 90.01 & 139.66 & 185.25 & 238.25 & 289.63 & 346.42 & 388.71 & 427.65 & 481.83 & 539.11 \\
$\Delta f$ (cents) & -37.82 & +8.22 & +2.34 & +31.17 & +70.95 & +78.83 & +80.74 & +43.67 & +50.43 & -53.41 \\
\multicolumn{11}{r}{\textit{RMSE: 52.48 cents}} \\
\hline \\
\multicolumn{11}{c}{\textbf{Reference alphorn}} \\
\hline
Notes & Do & Sol & Do & Mi & Sol & Sib & Do & Ré & Mi & Sol \\
\hline
$f_{\mathrm{note}}$ (Hz) & 92.00 & 139.00 & 185.00 & 234.00 & 278.00 & 331.00 & 371.00 & 417.00 & 468.00 & 556.00 \\
$f_{\mathrm{peak}}$ (Hz) & 94.11 & 146.83 & 196.97 & 246.98 & 296.51 & 343.50 & 388.71 & 436.02 & 487.77 & 537.15 \\
$\Delta f$ (cents) & +39.20 & +94.84 & +108.51 & +93.44 & +111.57 & +64.17 & +80.74 & +77.22 & +71.62 & -59.72 \\
\multicolumn{11}{r}{\textit{RMSE: 82.93 cents}} \\
\hline \\
\end{tabular}
\caption{Comparison of the reference alphorn geometry and the geometries obtained from the forward and backward models. For each case, the alphorn scale frequencies $f_{\mathrm{note}}$, the resonance frequencies $f_{\mathrm{peak}}$ extracted from FEM-computed impedance spectra, and the corresponding deviations $\Delta f$, expressed in cents relative to $f_{\mathrm{note}}$, are reported. The root mean squared error (RMSE) of the deviations is also provided for each configuration.}
\label{tab:final_comparison_models}
\end{table}

Overall, the agreement with the target frequencies is very good for the lowest notes and
gradually deteriorates as the frequency increases for both predicted geometries.
For the geometry predicted by the machine-learning models, the deviation remains
below $100$ cents (one semitone) for all target notes. This indicates that the produced pitches
would be perceived as close to the desired ones. However, some loss of accuracy is observed at higher frequencies.
Interestingly, the geometry predicted by the backward model yields resonance frequencies that are, except for the Re, closer to the target frequencies than those obtained from the measured alphorn geometry.
On the other hand, the geometry identified through the forward model consistently achieves deviations that are comparable to or smaller than those observed for the reference alphorn across all considered notes.
In particular, it exhibits remarkable performance for the first Sol and Do. Overall, this latter configuration provides the best global agreement with the target frequencies.

We remark that some of the observed deviations in cents are large enough to be perceptible, especially at higher frequencies. This is also the case for the impedance spectrum obtained from a real alphorn, which exhibits comparable deviations with respect to the official note frequencies. In practice, we measured the notes effectively on the real alphorn geometry and observed that they are significantly closer to the target frequencies than the corresponding impedance peaks
predicted by the numerical model. This may indicate that the ability to produce well-tuned notes is not solely determined by the exact alignment of impedance peaks with the target frequencies.

More specifically, impedance peaks describe preferred resonances of the linear acoustic resonator excited by a non-linear source (the lips).
In lip-reed instruments, the oscillation frequency can lock within a finite range around a resonance due to nonlinear interaction between the exciter and the acoustic impedance \cite{fletcher1978mode}. Since the
present model is linear and relies on a prescribed boundary condition at the mouthpiece, it cannot represent these locking ranges.
In particular, the mathematical model assumes a fixed normal velocity at the mouthpiece,
whereas a real player continuously adapts the airflow. It is therefore reasonable to assume
that, in such cases, a player can compensate for moderate discrepancies between
impedance peaks and desired notes by adjusting the airflow possibly at the cost of increased playing effort.

From this perspective, a geometry whose impedance peaks are closer to the target frequencies
may reduce the amount of active compensation required from the player. The geometry predicted by the machine-learning model could therefore be easier to play than the measured instrument. This interpretation remains limited as playability depends on complex player–instrument interactions that are not captured by the present acoustic model.
\section{Conclusions and Limitations}\label{sec:conclusions}
In this work, we investigated a computational framework combining finite element modeling, the reduced basis method and machine learning to study and design the acoustic behavior of the Swiss alphorn.
A central contribution of this study is the derivation of an affine, parametrized formulation on a fixed reference domain, enabling the efficient application of the reduced basis method. The use of frequency-localized reduced bases is crucial to control the approximation error over a broad frequency range.
The reduced-order simulations are used to generate a large dataset linking geometry to resonance frequencies, which in turn enabled the training of machine learning models for both forward and inverse acoustic mappings.
Overall, this work demonstrates that combining reduced-order modeling with machine
learning provides a viable and efficient approach for exploring and designing alphorn
geometries. While the proposed framework achieves substantial computational speedups, it
also shows the limitations of a linear model based on fixed boundary conditions when addressing
questions related to musical playability.

Future work could therefore focus on enriching the physical model by incorporating  additional effects such as viscothermal losses, wall vibrations, or more realistic boundary conditions at the mouthpiece. Another important direction would be to integrate player-in-the-loop information, either through experimental data or more advanced source models, to better capture the interaction between airflow control and acoustic response. Such
extensions would further improve the predictive power of the framework and strengthen its applicability to instrument design.

\appendix

\section{Appendix}
\subsection{Jacobian and affine decomposition computation}
In this appendix, we report the derivations of both the Jacobian of the cone transformation and the affine decomposition. 
\subsubsection{Jacobian of the transformation}\label{app:appendix_1}
The Jacobian matrix of $T_{\mu_c}$ has the explicit form :
\[
D T_{\mu_c}(x,y,z)
=
\begin{pmatrix}
\alpha(\mu_c) & 0 & 0 \\
\lambda'(x;\mu_c)\, y & \lambda(x;\mu_c) & 0 \\
\lambda'(x;\mu_c)\, z & 0 & \lambda(x;\mu_c)
\end{pmatrix},
\]
where:
\[
\lambda'(x;\mu_c)
=
\frac{\beta(\mu_c)\,\hat r_1 - \hat{\gamma} \, r_1}{D(x)^2}.
\]

\subsubsection{Affine decomposition of the Jacobian determinant}\label{app:appendix_2}
The determinant of the Jacobian is given by:
\[
\det\bigl(D T_{\mu_c}\bigr)
=
\alpha(\mu_c)\,\lambda(x;\mu_c)^2
\]

showing that the Jacobian determinant is a rational function of $x$,
with coefficients depending only on $\mu_c$.
Expanding
\[
\lambda(x;\mu_c)^2
=
\left( \frac{N(x;\mu_c)}{D(x)} \right)^2 \, ,
\]
yields
\[
\lambda(x;\mu_c)^2
=
\frac{r_1^2}{D(x)^2}
+
\frac{\,\beta(\mu_c)^2(x-\hat{x}_1)^2}{D(x)^2}
+
\frac{2r_1\beta(\mu_c)(x-\hat{x}_1)}{D(x)^2}.
\]

Consequently, the Jacobian determinant admits the affine decomposition
\[
\det\bigl(D T_{\mu_c}\bigr)
=
\alpha(\mu_c)\,
\left(
\frac{r_1^2}{D(x)^2}
+
\frac{\,\beta(\mu_c)^2(x-\hat{x}_1)^2}{D(x)^2}
+
\frac{2r_1\beta(\mu_c)(x-\hat{x}_1)}{D(x)^2}
\right).
\]

\subsubsection{Affine decomposition of the metric tensor}\label{app:appendix_3} 
We recall the metric tensor definition:
\[
G_{\mu_c}(\xi)
\;\coloneqq\;
\det\bigl(DT_{\mu_c}(\xi)\bigr)\,
DT_{\mu_c}(\xi)^{-1}
\bigl(DT_{\mu_c}(\xi)^{-1}\bigr)^{T}.
\]

A direct symbolic computation, using \texttt{SymPy} \cite{simpy2014simpy}  yields:
\[
G_{\mu_c}(x,y,z)
=
\frac{1}{\alpha(\mu_c)}
\begin{pmatrix}
\lambda(x;\mu_c)^2
&
- y \lambda(x;\mu_c) \lambda'(x;\mu_c)
&
- z \lambda(x;\mu_c) \lambda'(x;\mu_c)
\\[0.3em]
- y \lambda(x;\mu_c) \lambda'(x;\mu_c)
&
\alpha(\mu_c)^2 + y^2  \lambda'(x;\mu_c)^2
&
y z \lambda'(x;\mu_c)^2
\\[0.3em]
- z \lambda(x;\mu_c) \lambda'(x;\mu_c)
&
y z \lambda'(x;\mu_c)^2
&
\alpha(\mu_c)^2 + z^2 \lambda'(x;\mu_c)^2
\end{pmatrix}.
\]

This expression can be rewritten as :
\[
G_{\mu_c}
=
\alpha(\mu_c)
\begin{pmatrix}
0 & 0 & 0 \\
0 & 1 & 0 \\
0 & 0 & 1
\end{pmatrix}
+
\frac{\lambda(x;\mu_c)^2}{\alpha(\mu_c)}
\begin{pmatrix}
1 & 0 & 0 \\
0 & 0 & 0 \\
0 & 0 & 0
\end{pmatrix}
-
\frac{\lambda(x;\mu_c)\lambda'(x;\mu_c)}{\alpha(\mu_c)}
\begin{pmatrix}
0 & y & z \\
y & 0 & 0 \\
z & 0 & 0
\end{pmatrix}
+
\frac{\lambda'(x;\mu_c)^2}{\alpha(\mu_c)}
\begin{pmatrix}
0 & 0 & 0 \\
0 & y^2 & y z \\
0 & y z & z^2
\end{pmatrix}.
\]

Expanding the rational terms in \(x\) gives :
\[
\lambda(x;\mu_c)^2
=
\frac{r_1^2}{D(x)^2}
+
\frac{\beta(\mu_c)^2 (x-\hat{x}_1)^2}{D(x)^2}
+
\frac{2 r_1\,\beta(\mu_c)(x-\hat{x}_1)}{D(x)^2},
\]
\[
\lambda(x;\mu_c)\lambda'(x;\mu_c)
=
\frac{\beta(\mu_c)\hat r_1 r_1 - \hat{\gamma}r_1^{\,2}}{D(x)^3}
+
(\beta(\mu_c)^2 \hat r_1 - \beta(\mu_c)\hat{\gamma} r_1 )\frac{(x-\hat{x}_1)}{D(x)^3},
\]
\[
\lambda'(x;\mu_c)
=
\frac{\beta(\mu_c)\,\hat r_1 - \hat{\gamma} \, r_1}{D(x)^2}.
\]

We introduce the parameter-dependent coefficients:
\[
\omega_1(\mu_c) \coloneqq \alpha(\mu_c), \qquad
\omega_2(\mu_c) \coloneqq \frac{r_1^2}{\alpha(\mu_c)}, \qquad
\omega_3(\mu_c) \coloneqq \frac{2r_1\,\beta(\mu_c)}{\alpha(\mu_c)},
\]
\[
\omega_4(\mu_c) \coloneqq \frac{\beta(\mu_c)^2}{\alpha(\mu_c)}, \qquad
\omega_5(\mu_c) \coloneqq
- \frac{\beta(\mu_c)\hat r_1 r_1 - \hat{\gamma} r_1^{\,2}}{\alpha(\mu_c)},
\]
\[
\omega_6(\mu_c) \coloneqq
- \frac{\beta(\mu_c)^2 \hat r_1 - \beta(\mu_c)\hat{\gamma} r_1}{\alpha(\mu_c)}, \qquad
\omega_7(\mu_c) \coloneqq
\frac{\bigl(\beta(\mu_c)\hat r_1 - \hat{\gamma} r_1\bigr)^2}{\alpha(\mu_c)}.
\]

We define the seven parameter-independent matrices:
\[
A_1(x,y,z) \coloneqq
\begin{pmatrix}
0 & 0 & 0 \\
0 & 1 & 0 \\
0 & 0 & 1
\end{pmatrix},
\qquad
A_2(x,y,z) \coloneqq
\frac{1}{D(x)^2}
\begin{pmatrix}
1 & 0 & 0 \\
0 & 0 & 0 \\
0 & 0 & 0
\end{pmatrix},
\]
\[
A_3(x,y,z) \coloneqq
\frac{(x-\hat{x}_1)}{D(x)^2}
\begin{pmatrix}
1 & 0 & 0 \\
0 & 0 & 0 \\
0 & 0 & 0
\end{pmatrix},
\qquad
A_4(x,y,z) \coloneqq
\frac{(x-\hat{x}_1)^2}{D(x)^2}
\begin{pmatrix}
1 & 0 & 0 \\
0 & 0 & 0 \\
0 & 0 & 0
\end{pmatrix},
\]

\[
A_5(x,y,z) \coloneqq
\frac{1}{D(x)^3}
\begin{pmatrix}
0 & y & z \\
y & 0 & 0 \\
z & 0 & 0
\end{pmatrix},
\qquad
A_6(x,y,z) \coloneqq
\frac{x-\hat{x}_1}{D(x)^3}
\begin{pmatrix}
0 & y & z \\
y & 0 & 0 \\
z & 0 & 0
\end{pmatrix},
\]

\[
A_7(x,y,z) \coloneqq
\frac{1}{D(x)^4}
\begin{pmatrix}
0 & 0 & 0 \\
0 & y^2 & y z \\
0 & y z & z^2
\end{pmatrix}.
\]

With these definitions, the metric tensor admits the affine decomposition:
\[
G_{\mu_c}(x,y,z)
=
\sum_{i=1}^{7}
\omega_i(\mu_c)\, A_i(x,y,z).
\]

\subsection{Validation on benchmark problem}\label{apx:verification_on_test}
Before applying the model to the 3D alphorn geometry, we verify its correct implementation using a benchmark problem with a known analytical solution. The model problem, defined on the square $\Omega = [0, \pi]^2$, reads:

    \begin{align}
        \left\{
        \begin{aligned}
            &\Delta p + k^2 p = 0 \quad &&\text{in } \Omega \, , \\
            &\frac{\partial p}{\partial \mathbf{n}} = g_N \quad &&\text{for } x=0 \, ,\\
            &\frac{\partial p}{\partial \mathbf{n}} = 0 \quad &&\text{for } y=0, y=\pi \, ,\\
            &\frac{\partial p}{\partial \mathbf{n}} + \alpha p = 0 \quad &&\text{for } x=\pi \, .
        \end{aligned}
        \right.
        \label{eq:test}
    \end{align}

\begin{figure}[h!]
      \centering
        \includegraphics[width=0.45\textwidth]{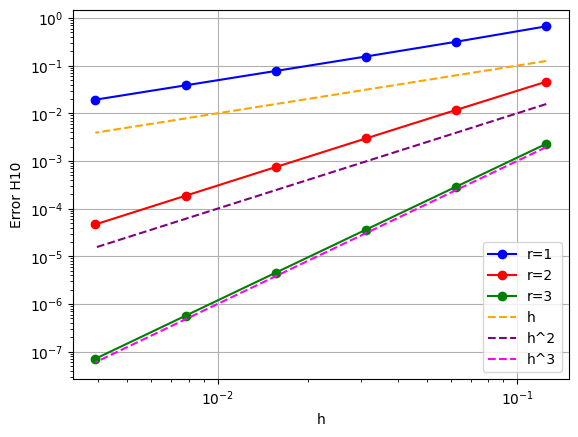} \quad
        \includegraphics[width=0.45\textwidth]{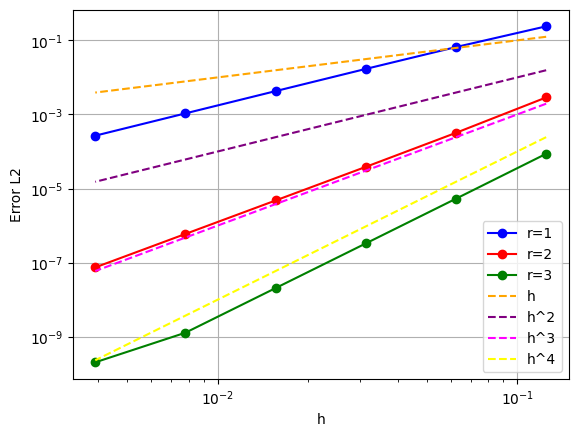}
        \caption{Error convergence for the validation test. h is the mesh size and r represents the polynomial degree of the FE basis. Left: Error convergence in $H^1$ norm against mesh refinements. Right: Error convergence in $L^2$ norm against mesh refinements.}
    \label{fig:error_convergence}
    \end{figure}

By imposing $k^2 = 2$, $\alpha = \sqrt{3}$ and $g_N(y) = \frac{\sqrt{3}}{2} \cos(y) + i \frac{\sqrt{3}}{2} \cos(y)$, the analytical solution of the system of equations (\ref{eq:test}) reads: $p_{ex}(x,y) = \cos(x+\frac{\pi}{3})\cos(y) + i\cos(x+\frac{\pi}{3})\cos(y)$.
This allows us to compute the numerical error and check the convergence of the method under mesh refinements. The results for error convergence in $L^2$ and $H^1$ norms are reported in Fig.\ \ref{fig:error_convergence}. As we can see, the implementation achieves the rates of convergence theoretically expected for smooth solutions, which for FEM are known to be:

    \begin{equation*}
        \| p_{ex} - p_h \|_{H^1} \lesssim h^r \, ,  \quad
        \| p_{ex} - p_h \|_{L^2} \lesssim h^{r+1} \, ,
    \end{equation*}

where $h = \displaystyle \max_{i} \{\text{diam}(K_i)\}$ indicates the mesh refinement.
    
Therefore, we conclude that the proposed FE implementation demonstrates consistent and reliable convergence under mesh refinement for the benchmark problem. These results validate the correctness and stability of the numerical scheme, indicating that the implementation provides a solid foundation for its use as the core solver in the Reduced Basis Method (RBM).

\subsection{Measurements of three real alphorns}\label{apx:real-alphorns-measurements}
\begin{table}[H]
\centering
\begin{tabular}{lccccccc}
\hline
Instrument & $L_1$ & $L_2$ & $L_3$ & $L_4$ & $R_2$ & $R_3$ & $R_4$ \\
\hline
Alphorn 1 & 544.96 & 844.87 & 754.84 & 444.94 & 13.50 & 28.50 & 44.00 \\
Alphorn 2 & 534.95 & 799.86 & 804.85 & 444.94 & 13.50 & 28.50 & 44.00 \\
Alphorn 3 & 539.92 & 844.91 & 754.84 & 444.94 & 16.00 & 28.50 & 44.00 \\
\hline \\
\end{tabular}
\caption{Geometrical parameters of the three alphorns (in mm).}
\label{tab:geom}
\end{table}

\subsection{Hyperparameter selection for the machine learning models}
The hyperparameters of both forward and backward models are selected using a combination of cross-validation and targeted experiments. The models are trained to minimize the mean squared error (MSE).

\subsubsection{Forward model}\label{apx:hyper-forward}
The learning rate is selected using a $4$-fold cross-validation procedure. The following values are tested:
\[
\{10^{-1},\,10^{-2},\,10^{-3},\,10^{-4},\,10^{-5}\},
\]
and the results show that a learning rate of $10^{-3}$ provides the best trade-off between convergence stability and accuracy.
The hidden dimension is explored in the range $\{32, 64, \dots, 256\}$. 
The results show a monotonic improvement with increasing width, without clear evidence of saturation in the explored range. 
The value of $256$ is retained as a trade-off between accuracy and computational cost.
The network depth is also selected using a validation study, where the following values are tested:
\[
\{2,\,3,\,4,\,5,\,6,\,7\}.
\]
The results show that shallow networks underfit the data, while deeper architectures beyond five layers do not provide significant improvements. The value of $5$ is therefore retained as a trade-off between accuracy and computational cost.

The final model consists of five hidden layers with $256$ neurons each and ReLU activation functions. Training is performed with mini-batches of size $256$, an $\ell_2$ regularization weight of $10^{-5}$, and $200$ training epochs.

\subsubsection{Backward model}\label{apx:hyper-backward}
For the backward model, the hidden dimension is selected using a $4$-fold cross-validation procedure. The following values are tested:
\[
\{128,\,256,\,512,\,1024,\,2048\}.
\]
The results indicate that increasing the width improves performance up to $1024$, whereas the validation error increases again for $2048$.
The final model consists of four hidden layers with $1024$ neurons each, ReLU activation functions, and a dropout rate of $0.1$ applied after each hidden layer. The output layer is designed so as to enforce the geometric constraints directly (as detailed in sec. \ref{sec:backward-model}). Training is performed using the Adam optimizer with a learning rate of $10^{-3}$ and an $\ell_2$ regularization weight of $10^{-3}$. As preprocessing, the input frequencies are first transformed using a base-10 logarithm and then standardized by subtracting the mean and dividing by the standard deviation computed on the training set. The same transformation is subsequently applied to the validation and test sets using these statistics.

\section*{Acknowledgments}
S. Deparis, F. Marcinn\`o and R. Tenderini have been supported by the
Swiss National Science Foundation under project "Data-driven approximation of haemodynamics by combined reduced order modeling and deep neural networks",  n. 200021-197021.

\bibliographystyle{plain}
\bibliography{refs}  

\end{document}